\documentclass[a4paper, 11pt]{article}
\usepackage{amsthm}
\usepackage{amsfonts}
\usepackage{latexsym}
\usepackage{epic}
\usepackage{graphics,multicol}
\usepackage{graphicx}%
\usepackage{epsf, epstopdf}
\usepackage{etex}
\usepackage{pinlabel}
\usepackage{ulem}

\usepackage[top=2.5cm,bottom=2.5cm,left=2.5cm,right=2.5cm]{geometry}
\usepackage{graphicx}
\usepackage{amsmath,amsthm,amssymb,lineno}

\usepackage{graphicx}           
\usepackage{color, ifpdf}  
\usepackage[dvipsnames]{xcolor}            
\usepackage{graphicx}
\usepackage{geometry}
\usepackage{float}
\usepackage{indentfirst}        
\usepackage{amsmath,amssymb,bm} 
\usepackage{cases}              
\usepackage{setspace}
\usepackage{graphicx,booktabs,multirow}
\usepackage{enumerate,mathdots}

\usepackage{latexsym}
\usepackage{graphicx,color, booktabs,multirow}
\usepackage{tikz,chemfig}
\usepackage{graphicx,booktabs,multirow,mathrsfs}
\usepackage{appendix}
\usepackage{yhmath}
\usepackage{setspace}

\newtheorem{theo}{Theorem}
\newtheorem{lemm}[theo]{Lemma}

\newtheorem{conj}[theo]{Conjecture}
\newtheorem{coro}[theo]{Corollary}
\newtheorem{prop}[theo]{Proposition}

\newtheorem{clm}{Claim}

\newcommand \Them[3]
{
	\begin{theo}[#1]
		\label{#2}	#3
	\end{theo}
}

\usetikzlibrary{backgrounds}
\usetikzlibrary{shapes,shapes.geometric,shapes.misc}
\pgfdeclarelayer{edgelayer}
\pgfdeclarelayer{nodelayer}
\pgfsetlayers{background,edgelayer,nodelayer,main}
\tikzstyle{none}=[inner sep=0mm]
\tikzstyle{every loop}=[]

\tikzstyle{dotted}=[dash pattern=on \pgflinewidth off 2pt]
\tikzstyle{dashed}=[dash pattern=on 3pt off 3pt]

\tikzstyle{new style 0}=[fill=black, draw=black, shape=circle]
\tikzstyle{red style 1}=[fill=red, draw=black, shape=circle]
\tikzstyle{blue style 2}=[fill=blue, draw=black, shape=circle]
\tikzstyle{white style 4}=[fill=white, draw=black, shape=circle]
\tikzstyle{bklack style 5}=[fill=black, draw=black, shape=rectangle]
\tikzstyle{red style 3}=[fill=red, draw=black, shape=rectangle]
\tikzstyle{yellow style 7}=[fill=yellow, draw=black, shape=rectangle]
\tikzstyle{new style 8}=[fill={rgb,255: red,0; green,132; blue,0}, draw={rgb,255: red,0; green,131; blue,0}, shape=circle]

\tikzstyle{new edge style 0}=[-]
\tikzstyle{new edge style 1}=[-, draw=red]
\tikzstyle{new edge style 2}=[-, draw=blue]
\tikzstyle{new edge style 3}=[-, draw={rgb,255: red,0; green,156; blue,0}]

\allowdisplaybreaks

\numberwithin{equation}{section}

\newcounter{countcase}

\newcounter{countclaim}
\def\inclaim{\addtocounter{countclaim}{1}
{\noindent {\bf Claim \thecountclaim}: }}

\def \proof {\noindent {\it Proof}. }
\newcommand{\claimend}{{\hfill$\natural $}}

\newcommand{\proofend}{
{\hfill$\Box$}
\setcounter{countclaim} {0}
\setcounter{countcase} {0}
}

\def \N {{\mathbb {N}}}
\def\P {{\mathbb {P}}}
\def \setm {{\cal M}}
\def \sets {{\cal S}}
\def\setd {{\cal D}}
\def \iff  {{if and only if }}

\newcommand \equ[2]
{
\begin{equation}\label{#1}
#2
\end{equation}
}

\newcommand \eqn[2]
{
\begin{eqnarray}\label{#1}
#2
\end{eqnarray}
}

\newcommand \aln[2]
{
\begin{align}\label{#1}
#2
\end{align}
}

\newcommand \upo[1]
{
	(#1\uparrow )
}

\newcommand \downo[1]
{
	(#1\downarrow )
}

\begin{document}
		\baselineskip 0.6 cm
	
	\title{
		The absolute values of
	the perfect matching derangement graph's eigenvalues
	 almost follow the lexicographic order of partitions
	}
	
	\author{
		Meiqiao Zhang\thanks{Corresponding Author. Email: nie21.zm@e.ntu.edu.sg and 
			meiqiaozhang95@163.com.},
		Fengming Dong\thanks{
			Email: fengming.dong@nie.edu.sg and donggraph@163.com.}
		\\
		National Institute of Education,
		Nanyang Technological University, Singapore
	}
	\date{}
	
	\maketitle{}

	\begin{abstract}
		In 2013, Ku and Wong showed 
	that for any partitions $\mu$ and $\mu'$ of a positive integer $n$
		with the same first part $u$
		and the lexicographic order $\mu\triangleleft \mu'$, 
		the eigenvalues $\xi_{\mu}$ and $\xi_{\mu'}$ of 	
		the derangement graph $\Gamma_n$ have the property 
		 $|\xi_{\mu}|\le |\xi_{\mu'}|$,
		 where the equality holds \iff 
		$u=3$ and all other parts are less than $3$.
		In this article, we obtain an analogous conclusion 
		on the eigenvalues of the perfect matching derangement graph $\setm_{2n}$ of $K_{2n}$
		by  finding 
		a new recurrence formula 
		for the eigenvalues of  $\setm_{2n}$. 
		\end{abstract} 

	\noindent {\bf Keywords:}
	Cayley graph, the derangement graph,
the perfect matching derangement graph, eigenvalue
	
		\smallskip
	\noindent {\bf Mathematics Subject Classification: 05A17, 05C50}

	\section{ Introduction
		\label{secintro}}
	
	Let $\N$ denote the set of positive integers. For any $n\in \N$, let $\sets_n$ be the symmetric group on $[n]$, where 
	$[n]:=\{1,2,\dots,n\}$, and let $\setd_n$ be the set of derangements in $\sets_n$, 
	where $\pi\in \sets_n$ is called a 
	{\it derangement} if $\pi(i)\ne i$ for all $i\in [n]$. 
	Any two members $\pi_1,\pi_2$ of $ \sets_n$ are said to be 
	a derangement of each other
	if $\pi_1(i)\ne \pi_2(i)$ for all 
	$i\in [n]$.
	The {\it derangement graph} $\Gamma_n$ is defined to be the graph
	with vertex set $\sets_n$ 
	in which any two vertices 
	$\pi_1$ and $\pi_2$ are adjacent 
	\iff $\pi_1$ is a derangement of $\pi_2$.
	 Obviously, $\Gamma_n$ is the Cayley graph $Cay(\sets_n,\setd_n)$.
	 
	 For any $n\in \N$,  a \textit{partition} $\mu$ of $n$, denoted by $\mu\vdash n$, 
	 is a non-increasing sequence of positive integers $(\mu_1,\dots,\mu_r)$
	 such that 
	 $\mu_1+\cdots+\mu_r=n$.
	 The {\it length of $\mu$}, denoted by 
	 $\ell(\mu)$, is defined to be $r$, 
	 the {\it size of $\mu$}, denoted by 
	 $|\mu|$, is defined to be $n$,
	 and $\mu_i$ is called the $i$-th part of $\mu$.

	 Since $\setd_n$ is closed under conjugation, 
	the Cayley graph $\Gamma_n=Cay(\sets_n,\setd_n)$  is normal.
	Based on the fact that 
	 the conjugacy classes of $\sets_n$ and the irreducible characters of
	 $\sets_n$ are both indexed by 
	 partitions  of $n$,
	 the eigenvalues of  
	 $\Gamma_n$ can be 
	 denoted by $\xi_{\mu}$,
	 where $\mu\vdash n$
	 (see \cite{koh2023,ku2010,ku2013, Rente2007}).

	 For any  $\mu=(\mu_1,\dots,\mu_r)
	 \vdash n$, define
	 $
	 \mu\setminus \mu_r:=
	 (\mu_1,\dots,\mu_{r-1})\vdash (n-\mu_r)
	 $
	 and for $1\le k\le \mu_r$, define
	 $
	 \mu-\hat{k} 
	 :=(\mu_1-k,\dots,\mu_{r}-k)
	 \vdash (n-kr).
	 $
	 In particular,
	 when $k=\mu_r$, 
	 the partition $\mu-\hat{k}$ might have many parts equal to $0$.
	  Thus, we assume that 
	 $
	 	 (\mu_1,\mu_2,\dots,\mu_i,0,\dots,0)
	 	 =(\mu_1,\mu_2,\dots,\mu_i)
	 	 $
	 	 and
	 $(0,0,\dots,0)=(0)$.	
	 Renteln~\cite{Rente2007} found 
	 a recurrence formula for the eigenvalues $\xi_{\mu}$ of $\Gamma_n$,
	 where $\xi_{(0)}:=1$ by convention.

	 \Them{Renteln~\cite{Rente2007}}
	 {the-1} 
	 {
	 	For any $\mu =(\mu_1,\dots,\mu_r)\vdash n$ with $r\ge 2$, 
	 	$$
	 	\xi_{\mu}=(-1)^{r-1}(\mu_1+r-1) \xi_{\mu-\hat{1}} +(-1)^{\mu_1+r-1}
	 	\xi_{(\mu_2-1,\dots, \mu_r-1)}.
	 	$$ 
	 }
Ku and Wong \cite{ku2013} obtained a similar recurrence formula for the eigenvalues $\xi_{\mu}$ of $\Gamma_n$. 
	 
	  \Them{Ku and Wong~\cite{ku2013}}
	 {the0} 
	 {
	 For any $\mu =(\mu_1,\dots,\mu_r)\vdash n$ with $r\ge 2$, 
	 $$
	 \xi_{\mu}=(-1)^{r-1}\mu_r \xi_{\mu-\hat{1}} +(-1)^{\mu_r}\xi_{\mu\setminus \mu_r-\hat{1}}.
	 $$ 
 }
	 
	 For $u, n\in \N$ with $0<u\le n$, 
	 let $\P_{n,u}$ be the set of partitions $\mu\vdash n$ such that the first part of $\mu$ is equal to $u$.
	 For partitions $\mu=(\mu_1,\dots,\mu_r)$
	 and $\mu'=(\mu'_1,\dots,\mu'_s)$,
	 we say $\mu$ is {\it dominated} by $\mu'$,
	 written as $\mu\trianglelefteq \mu'$, 
	 if $\mu_1+\cdots+\mu_k\le \mu'_1+\cdots+\mu'_k$
	for all $k\in [r]$, where $\mu_k'=0$ if $k>s$.
	 Write $\mu\triangleleft \mu'$ if 
	 $\mu\ne \mu'$ and 
	 $\mu\trianglelefteq \mu'$.

	 By applying Theorem~\ref{the0}, Ku and Wong further proved 
	 the following result on the 
	 absolute values
	 of eigenvalues of 
	 $\Gamma_n$.
	 For any $n\ge 3$, let $\P^*_{n,3}$ be the set of partitions $\mu=(\mu_1,\dots,\mu_r)\in \P_{n,3}$ 
	 with $1\le \mu_2\le 2$ whenever $r\ge 2$.

	 \Them{Ku and Wong~\cite{ku2013}}
	 {the1} 
	 {
	 	For partitions $\mu,\mu'\in \P_{n,u}$, 
	 if $\mu \triangleleft \mu'$, then 
	$
	 |\xi_{\mu}| \le  |\xi_{\mu'} |.
	 $
	 Furthermore, the inequality is tight if and only if  $u=3$  and 
	 $\mu,\mu'\in \P^*_{n,3}$.
	}

 Theorem~\ref{the1} actually confirmed 
 the following conjecture.

\begin{conj}[Ku and Wales \cite{ku2010}]\label{con1}
For any $\mu=(\mu_1,\dots,\mu_r)\vdash n$, let $\mu^*$ be the largest partition in lexicographic order among all the partitions in $\P_{n,\mu_1}$.
Then
$$
|\xi_{(\mu_1,1^{n-\mu_1})}|\le |\xi_{\mu}|\le |\xi_{\mu^*}|,
$$
where $(\mu_1,1^{n-\mu_1})$ is the 
partition in $\P_{n,\mu_1}$ 
with the $i$-th part equal to $1$ 
for all $i: 2\le i\le n-\mu_1+1$. 
\end{conj}

In this article, we will show that the eigenvalues of the perfect matching derangement graph
have a similar property 
as Theorem~\ref{the1}.

For any $n\in \N$, let $X_{2n}$ denote the set of perfect matchings of the complete graph $K_{2n}$.
It can be verified that $|X_{2n}|=(2n-1)!!$, 
where $(2n-1)!!=1\cdot 3 \cdots (2n-1)$.
The {\it perfect matching derangement graph} with respect to $K_{2n}$,
denoted by $\setm_{2n}$, 
is the graph with vertex set $X_{2n}$
such that any two vertices $M_1$
and $M_2$ are adjacent \iff 
$M_1\cap M_2=\emptyset$, 
i.e., $M_1$ and $M_2$ are derangements with respect to each other.
Then due to symmetry, $\setm_{2n}$ is a regular graph, whose degree $d_n$ can be shown to be equal to $\sum\limits_{i=0}^{n-1}(-1)^i\binom{n}{i} (2n-2i-1)!!$ by the principle of inclusion-exclusion. Meanwhile, a recurrence relation for calculating $d_n$ emerges as
\aln{the5-e2}
{
d_n=2(n-1)(d_{n-1}+d_{n-2}).
}

It is known that the eigenvalues of 
$\setm_{2n}$ can also be 
indexed by partitions $\lambda$ of $n$
(see \cite{god2017,koh2023,lind2017,Rente2022}), and are denoted by 
$\eta_{\lambda}$.
When $\ell(\lambda)=1$ (i.e., $\lambda=(n)$), $\eta_{\lambda}=d_n$. For the case when $\ell(\lambda)\ge 2$,
the authors of~\cite{koh2023} found 
the following recurrence formula 
for $\eta_{\lambda}$, where $\eta_{(0)}:=1$ by convention.

\Them{Koh, Ku and Wong~\cite{koh2023}}
{the2}
{
The eigenvalues of the perfect matching derangement graph satisfy the following recurrence relation:
\equ{the2-e1}
{
(-1)^{\lambda_r}\eta_{\lambda}
=\eta_{\lambda \setminus \lambda_r}
+\sum_{j=1}^{\lambda_r}
(-1)^{jr} {\lambda_r\choose j}(2j-1)!!
\eta_{\lambda \setminus \lambda_r-\hat{j}},
}
where $\lambda=(\lambda_1,\dots, \lambda_r)$ and $r\ge2$. 
}

By applying the recurrence formula in Theorem~\ref{the2}, 
the authors in \cite{koh2023}  
proved 
the following alternating sign property
of eigenvalues of $\setm_{2n}$,
which was conjectured in \cite{god2015}
and \cite{lind2017}. 

\Them{Koh, Ku and Wong~\cite{koh2023}} 
{the3}
{
For $n\in\N$ with $n\ge 2$,
the perfect matching derangement graph $\setm_{2n}$
satisfies the alternating sign property, i.e., $(-1)^{n-\lambda_1}\eta_{\lambda}>0$ holds
for any $\lambda =(\lambda_1,\dots, \lambda_r)\vdash n$.
}

Note that Theorem~\ref{the3} was independently proved by  Rentel~\cite{Rente2022}  using different techniques.

In this article, we will first apply (\ref{the2-e1}) to obtain a new 
recurrence formula for $\eta_{\lambda}$. 

\Them{}{the4}
{
Let $\lambda=
(\lambda_1,\dots,\lambda_{s})\vdash n$,
where $s\ge 2$.
For any $2\le i\le s$, 
if either $i=s$ or $\lambda_{i}>\lambda_{i+1}$,
then 
the following recurrence relation 
holds:  
\aln{the4-e1}
{
	\eta_{\lambda}
	=-\eta_{\lambda'}
	+(-1)^{s+1} (2\lambda_i+s-i-1)
	\eta_{\lambda-\hat{1}}
	+(-1)^{s+1} (2\lambda_i+s-i-2)
	\eta_{\lambda'-\hat{1}},
}
where 
$\lambda'$ is the partition 
$(\lambda_1,\dots,\lambda_{i-1},
\lambda_{i}-1, \lambda_{i+1},
\dots,\lambda_{s})$ of $n-1$. 
}

Further, applying Theorem~\ref{the4}, we will 
obtain the following conclusion 
on the absolute values
of eigenvalues of 
$\setm_{2n}$ 
which is analogous to Theorem~\ref{the1}. 

\Them{}{the6}
{
For partitions $\lambda,\lambda'\in 
	\P_{n,u}$,
	if $\lambda \triangleleft \lambda'$,
	then
$|\eta_{\lambda}|\le|\eta_{\lambda'}|$,
 where the equality holds if and only if $u=3$ and $\lambda,\lambda'\in \P^*_{n,3}$.
}

As a special case, the next conclusion follows directly from Theorem~\ref{the6}.

\begin{coro}\label{cor1}
For any $\lambda=(\lambda_1,\dots,\lambda_r)\vdash n$, let $\lambda^*$ be the largest partition in lexicographic order among all the partitions in $\P_{n,\lambda_1}$.
	Then
	$
	|\eta_{(\lambda_1,1^{n-\lambda_1})}|\le |\eta_{\lambda}|\le |\eta_{\lambda^*}|.
	$

\end{coro}

Theorems~\ref{the4} and~\ref{the6}
will be proved in the 
following sections. 
	
\section{Proof of Theorem~\ref{the4}
		\label{sectwo}}

Recall that for  any $n\in \N$,  $d_n$ 
is the degree of the regular graph $\setm_{2n}$.  Specially, define $d_0=1$. 
Thus, the sequence 
$(d_0,d_1,d_2,\dots)$ 
can be determined by the recurrence relation (\ref{the5-e2}) with the initial conditions 
 $d_0=1$ and  $d_1=0$.

We now define a function $f$ on 
partitions of positive integers. 
For any $\lambda=(\lambda_1,\dots,\lambda_r)\vdash n$, 
define $f(\lambda)
=(-1)^{n-\lambda_1}\eta_{\lambda}$ 
and write $f(\lambda)$ as $f(\lambda_1,\dots,\lambda_r)$.
Thus, if $\ell(\lambda)=1$, i.e., $\lambda=(n)$, then 
$f(\lambda)=d_{n}$.
If $r\ge 2$, 
by Theorem~\ref{the2}, 
the following recurrence relation holds:
\aln{def1}
{
	f(\lambda)=f(\lambda\setminus \lambda_r)+\sum_{k=1}^{\lambda_r}\binom{\lambda_r}{k}(2k-1)!!
	f(\lambda\setminus \lambda_r-\hat{k}),
}
where 
$
f(\lambda_1,\lambda_2,\dots,\lambda_i,0,\dots,0)
:=f(\lambda_1,\lambda_2,\dots,\lambda_i)
$
for any partition $(\lambda_1,\lambda_2,\dots,\lambda_i)$ and $f(0):=d_0=1$.
Then it is clear that $f(\lambda)\ge 0$, and the equality holds if and only if $\lambda=(1)$.

For $r\ge 2$ and $2\le i\le r$, let 
$\P_n(r,i)$ denote the set of partitions 
$\lambda=(\lambda_1,\dots,\lambda_r)
\vdash n$ with 
$\lambda_{i-1}>\lambda_i$.
For any $\lambda=(\lambda_1,\dots,\lambda_r)
\in \P_n(r,i)$,
let $\lambda\upo{i}$ denote the 
partition obtained from $\lambda$ 
by replacing $\lambda_i$ by $\lambda_i+1$.
Clearly, $\lambda\upo{i}\vdash n+1$.
For any $\mu=(\mu_1,\dots,\mu_r)\vdash n$ and $2\le i\le r$,
if either $i=r$ or 
$\mu\in \P_n(r,i+1)$, let
$\mu\downo{i}$ denote the 
partition obtained from $\mu$ by
replacing $\mu_i$
by $\mu_i-1$.
Thus,  $\mu\downo{i}\vdash n-1$.

We first give a useful lemma for later calculations based on the recurrence relation (\ref{def1}).

\begin{lemm}\label{pp3}
For any  $\mu=(\mu_1,\dots,\mu_s)\vdash n$ with $s\ge 2$,
$$
\sum_{k=1}^{\mu_s}\binom{\mu_s}{k}(2k+1)!!f(\mu \setminus \mu_s-\hat{k})
=(2\mu_s+1)f(\mu)-2\mu_s
f(\mu\downo{s}) 
-f(\mu\setminus \mu_{s}).
$$
\end{lemm}

\proof If $\mu_s=1$, then 
by (\ref{def1}), 
$$
f(\mu)=f(\mu_1,\dots,\mu_{s-1},1)=f(\mu \setminus \mu_s)+f(\mu \setminus \mu_s -\hat{1}),
$$
and the result follows directly.

In the following, assume $\mu_s\ge 2$.
Note that 
$$
(2k+1)=(2\mu_s+1)+(2k-2\mu_s)
\quad  \mbox{and} \quad 
\binom{\mu_s}{k}(\mu_s-k)=\mu_s\binom{\mu_s-1}{k}.
$$
 Then
\eqn{eq1-2}
{
& & \sum_{k=1}^{\mu_s}\binom{\mu_s}{k}(2k+1)!!f(\mu\setminus \mu_s-\hat{k})
\nonumber\\
&=&
(2\mu_s+1)\sum_{k=1}^{\mu_s}\binom{\mu_s}{k}(2k-1)!!
f(\mu\setminus \mu_s-\hat{k})
+
\sum_{k=1}^{\mu_s-1}\binom{\mu_s}{k}(2k-2\mu_s)(2k-1)!!
f(\mu\setminus \mu_s-\hat{k})
\nonumber \\ 
&=&
(2\mu_s+1)\sum_{k=1}^{\mu_s}\binom{\mu_s}{k}(2k-1)!!
f(\mu\setminus \mu_s-\hat{k})
-
2\mu_s \sum_{k=1}^{\mu_s-1}\binom{\mu_s-1}{k}(2k-1)!!
f(\mu\setminus \mu_s-\hat{k})
\nonumber\\
&=&
(2\mu_s+1)\big (f(\mu)
-f(\mu\setminus \mu_s)\big )
-
2\mu_s\big (f(\mu\downo{s})
-f(\mu\setminus \mu_s)\big )
\nonumber\\
&=&
(2\mu_s+1)f(\mu)-
2\mu_sf(\mu\downo{s})-f(\mu\setminus \mu_s),
}
where
the second last expression follows from  (\ref{def1}).
The result holds.
\proofend

In the following, we focus on finding 
an explicit expression of $f(\mu\upo{i})-
f(\mu)$ for any partition $\mu\in \P_n(r,i)$.
We first deal with the simplest case 
$i=r$.

\begin{lemm}\label{pp2}
For $\mu=(\mu_1,\dots,\mu_s)\in\P_n(s,s)$,
\aln{eq2-1}
{
f(\mu\upo{s})-
f(\mu)=(2\mu_s+1)f(\mu\upo{s}-\hat{1})-2\mu_sf(\mu-\hat{1}).
}
\end{lemm}

\proof By the assumption $\mu\in\P_n(s,s)$,  $\ell(\mu\upo{s})=s\ge 2$ and $\mu_{s-1}>\mu_s\ge 1$.
By (\ref{def1}), 
\eqn{eq1-1}
{
f(\mu\upo{s})-f(\mu)
&=& \sum_{k=1}^{\mu_s+1} \binom{\mu_s+1}{k}(2k-1)!!
f(\mu\setminus \mu_s-\hat{k})
-
\sum_{k=1}^{\mu_s} \binom{\mu_s}{k}(2k-1)!!
f(\mu\setminus \mu_s-\hat{k})
\nonumber\\
&=&(2\mu_s+1)!!f
(\mu\setminus \mu_s-\widehat{(\mu_s+1)})+ \sum_{k=1}^{\mu_s} \binom{\mu_s}{k-1}(2k-1)!!
f(\mu\setminus \mu_s-\hat{k})
\nonumber\\
&=&(2\mu_s+1)!!f
(\mu\setminus \mu_s-\widehat{(\mu_s+1)})+ \sum_{j=0}^{\mu_s-1} \binom{\mu_s}{j}(2j+1)!!
f(\mu\setminus \mu_s-\widehat{(j+1)})
\nonumber\\
&=& f(\mu\setminus \mu_s- \hat{1})+\sum_{j=1}^{\mu_s} \binom{\mu_s}{j}(2j+1)!!
f(\mu\setminus \mu_s-\widehat{(j+1)})
\nonumber\\
&=& f(\mu\setminus \mu_s- \hat{1})+\sum_{j=1}^{\mu_s} \binom{\mu_s}{j}(2j+1)!!
f((\mu\upo{s}-\hat{1})\setminus \mu_s-\hat{j}).
}

Then by Lemma~\ref{pp3}, we have
\eqn{eq1-1-1}
{
f(\mu\upo{s})-f(\mu)
&=& f(\mu\setminus \mu_s-\hat{1})+
(2\mu_s+1)
f(\mu\upo{s} -\hat{1})
-2\mu_sf(\mu-\hat{1})-f(\mu\setminus \mu_s-\hat{1})
\nonumber\\
&=&
(2\mu_s+1)f(\mu\upo{s}-\hat{1})
-2\mu_sf(\mu-\hat{1}).
}
The result holds.
\proofend

Now, to handle the general cases, we prepare the following Proposition~\ref{pp4}, which is a key step towards proving Theorem~\ref{the4}.

\begin{prop}\label{pp4}
For any  $\mu=(\mu_1,\dots,\mu_{s})\in\P_n(s,i)$, where $2\le i\le s$,
\equ{eq4-1}
{
f(\mu\upo{i})-f(\mu)
	=(2\mu_i+s-i+1)
	f(\mu\upo{i}-\hat{1}) -(2\mu_i+s-i)
	f(\mu-\hat{1}).
}
\end{prop}

\proof
We shall prove the result by induction on $s-i$. 

Note that the case when $s-i=0$ is proven in Lemma~\ref{pp2}. Assume the result holds whenever 
$s-i<t$, where $0<t\le s-2$.
Now consider the case $s-i=t$.
Thus, $2\le i<s$.

We first deal with the case when $\mu_s=1$. 
By (\ref{def1}),
\eqn{eq4-1-0}
{f(\mu\upo{i})-f(\mu)
&=& \big (f(\mu\upo{i}\setminus \mu_s)
+f(\mu\upo{i}\setminus \mu_s-\hat{1})\big )
-\big(  f(\mu\setminus \mu_s)
+f(\mu\setminus \mu_s-\hat{1})  \big) 
\nonumber\\
&=&\big ( f(\mu\upo{i}\setminus \mu_s)-f(\mu\setminus \mu_s)\big)
+\big ( f(\mu\upo{i}\setminus \mu_s-\hat{1})
-f(\mu\setminus \mu_s-\hat{1}) \big)
\nonumber\\
&=&\big [ f(\mu\upo{i}\setminus \mu_s)-f(\mu\setminus \mu_s)\big]
+\big [ f(\mu\upo{i}-\hat{1})
-f(\mu-\hat{1}) \big ],
}
where the last expression follows from the fact that $\mu_s=1$.

Note that $\mu\upo{i}\setminus \mu_s
=(\mu \setminus \mu_s)\upo{i}$ 
as $i<s$.
By induction, we have
\eqn{eq4-1-1}
{
&&
f(\mu\upo{i}\setminus \mu_s)-f(\mu\setminus \mu_s)
\nonumber\\
&=& f((\mu \setminus \mu_s)\upo{i})-f(\mu\setminus \mu_s)
\nonumber\\
&=&(2\mu_i+(s-1-i)+1)
f((\mu \setminus \mu_s)\upo{i}
-\hat{1})
-(2\mu_i+(s-1-i))f(\mu\setminus \mu_s-\hat{1})
\nonumber\\
&=&(2\mu_i+(s-1-i)+1)
f(\mu\upo{i}-\hat{1})
-(2\mu_i+(s-1-i))f(\mu-\hat{1}),
}
where the last expression also follows from the fact that $\mu_s=1$.

Hence,  when $\mu_s=1$, 
the result follows from 
(\ref{eq4-1-0}) and (\ref{eq4-1-1}).

In the following, assume that 
$\mu_s\ge 2$.
By (\ref{def1}),  we have 
\eqn{eq4-2}
{
&&f(\mu\upo{i})-f(\mu)
\nonumber\\
&=&f(\mu\upo{i}\setminus \mu_s)
-f(\mu\setminus \mu_s)
\nonumber\\
&& +\sum_{k=1}^{\mu_s}
\left [ \left(\binom{\mu_s-1}{k}+\binom{\mu_s-1}{k-1}\right)(2k-1)!!
\big ( f(\mu\upo{i}\setminus \mu_s-\hat{k})
-f(\mu\setminus \mu_s-\hat{k}) \big)
\right ]
\nonumber \\
&=&f((\mu\setminus \mu_s)\upo{i})
-f(\mu\setminus \mu_s)
\nonumber\\
&& +\sum_{k=1}^{\mu_s-1}
 \binom{\mu_s-1}{k}(2k-1)!!
\big( f((\mu\setminus \mu_s)\upo{i}-\hat{k})
-f(\mu\setminus \mu_s-\hat{k}) 
\big )
\nonumber\\
&& +\sum_{k=1}^{\mu_s}
\binom{\mu_s-1}{k-1}(2k-1)!!
\big( f(\mu\upo{i}\setminus \mu_s-\hat{k})
-f(\mu\setminus \mu_s-\hat{k}) \big ),
}
where the last equality follows from the fact that $i<s$.

By induction, we have
\eqn{eq4-2-5}
{
& &f((\mu\setminus \mu_s)\upo{i})
-
f(\mu\setminus \mu_s)
\nonumber \\
&=&(2\mu_i+s-1-i+1)
f((\mu\setminus \mu_s)\upo{i}-\hat{1})
-(2\mu_i+s-1-i)f(\mu\setminus \mu_s-\hat{1})
\nonumber \\
&=&(2\mu_i+s-i)f(\mu\upo{i}\setminus \mu_s-\hat{1})
-(2\mu_i+s-i-1)f(\mu\setminus \mu_s-\hat{1}),
}
and similarly, for any $k$ with 
$1\le k\le \mu_s-1$, 
\eqn{eq4-2-5-0}
{
f((\mu\setminus \mu_s)\upo{i}-\hat{k})
-f(\mu\setminus \mu_s-\hat{k})
&=&(2(\mu_i-k)+s-i)f(\mu\upo{i}\setminus \mu_s-\widehat{(k+1)})
\nonumber \\
& &-(2(\mu_i-k)+s-i-1)f(\mu\setminus \mu_s-\widehat{(k+1)}),
}
implying that 
\eqn{eq4-2-6}
{
& &\sum_{k=1}^{\mu_s-1}
\binom{\mu_s-1}{k}(2k-1)!!
\big (f((\mu\setminus \mu_s)\upo{i}-\hat{k})
-f(\mu\setminus \mu_s-\hat{k}) \big )
\nonumber \\
&=&\sum_{k=1}^{\mu_s-1}
\binom{\mu_s-1}{k}(2k-1)!!
\big(
(2\mu_i+s-i)f(\mu\upo{i}\setminus \mu_s-\widehat{(k+1)})
-(2\mu_i+s-i-1)f(\mu\setminus \mu_s-\widehat{(k+1)})
\big)
\nonumber \\
& &+\sum_{k=1}^{\mu_s-1}
\binom{\mu_s-1}{k}(2k-1)!!
(-2k)\big(
f(\mu\upo{i}\setminus \mu_s-\widehat{(k+1)})
-f(\mu\setminus \mu_s-\widehat{(k+1)})
\big).
}

By (\ref{eq4-2}),  (\ref{eq4-2-5})
and  (\ref{eq4-2-6}), we have 
\equ{eq4-2-7}
{
f(\mu\upo{i})-f(\mu)=Q_1+Q_2,
}
where
\eqn{eq4-2-8}
{
Q_1&=&(2\mu_i+s-i)
\left (
f(\mu\upo{i}\setminus \mu_s-\hat{1})
+\sum_{k=1}^{\mu_s-1}
\binom{\mu_s-1}{k}(2k-1)!!
f(\mu\upo{i}\setminus \mu_s-\widehat{(k+1)})
\right )
\nonumber \\
& &-(2\mu_i+s-i-1)
\left (
f(\mu\setminus \mu_s-\hat{1})
+\sum_{k=1}^{\mu_s-1}
\binom{\mu_s-1}{k}(2k-1)!!
f(\mu\setminus \mu_s-\widehat{(k+1)})
\right )
\nonumber\\
&=&(2\mu_i+s-i)
\left (
f((\mu\upo{i}-\hat{1})\setminus (\mu_s-1))
+\sum_{k=1}^{\mu_s-1}
\binom{\mu_s-1}{k}(2k-1)!!
f((\mu\upo{i}-\hat{1})\setminus (\mu_s-1)-\hat{k})
\right )
\nonumber \\
& &-(2\mu_i+s-i-1)
\left (
f((\mu-\hat{1})\setminus (\mu_s-1))
+\sum_{k=1}^{\mu_s-1}
\binom{\mu_s-1}{k}(2k-1)!!
f((\mu-\hat{1})\setminus (\mu_s-1)-\hat{k})
\right )
\nonumber \\
&=&(2\mu_i+s-i)f(\mu\upo{i}-\hat{1})
-(2\mu_i+s-i-1)f(\mu-\hat{1}),
}
where the last expression follows from 
(\ref{def1}),
and 
\eqn{eq4-2-9}
{
Q_2&=&
\sum_{k=1}^{\mu_s-1}
\binom{\mu_s-1}{k}(2k-1)!!
(-2k)\big(
f(\mu\upo{i}\setminus \mu_s-\widehat{(k+1)})
-f(\mu\setminus \mu_s-\widehat{(k+1)})
\big)
\nonumber\\
&& +\sum_{k=1}^{\mu_s}
\binom{\mu_s-1}{k-1}(2k-1)!!
\big( f(\mu\upo{i}\setminus \mu_s-\hat{k})
-f(\mu\setminus \mu_s-\hat{k}) \big )
\nonumber \\ &=&
\sum_{k=1}^{\mu_s-1}
\binom{\mu_s-1}{k}(2k-1)!!
(-2k)\big(
f(\mu\upo{i}\setminus \mu_s-\widehat{(k+1)})
-f(\mu\setminus \mu_s-\widehat{(k+1)})
\big)
\nonumber\\
&& +\sum_{j=0}^{\mu_s-1}
\binom{\mu_s-1}{j}(2j+1)!!
\big( f(\mu\upo{i}\setminus \mu_s-\widehat{(j+1)})
-f(\mu\setminus \mu_s-\widehat{(j+1)}) \big )
\nonumber\\
&=&f(\mu\upo{i}\setminus \mu_s-\hat{1})- 
f(\mu\setminus \mu_s-\hat{1})
\nonumber \\
& &+
\sum_{k=1}^{\mu_s-1}
\binom{\mu_s-1}{k}(2k-1)!!\big(
f(\mu\upo{i}\setminus \mu_s-\widehat{(k+1)})
-f(\mu\setminus \mu_s-\widehat{(k+1)})
\big)
\nonumber \\
&=&f(\mu\upo{i}-\hat{1})-f(\mu-\hat{1}),
}
where the last expression follows from 
(\ref{def1}).

By (\ref{eq4-2-7}), (\ref{eq4-2-8}) and (\ref{eq4-2-9}),
 (\ref{eq4-1}) follows for the case $i\ge 2$ and $\mu_s\ge 2$.
 \proofend 
 
\textbf{Remark.}
 Proposition~\ref{pp4} does not hold for the case when $i=1$.
 This is because the initial step $f(\mu_1+1)-f(\mu_1)= (2\mu_1+1)f(\mu_1)-2\mu_1 f(\mu_1-1)$ fails due to (\ref{the5-e2}).

 Now we conclude this section by proving Theorem~\ref{the4}.
 
 \noindent \textit{Proof of Theorem~\ref{the4}.}
By 
 the definition,  $\eta_{\lambda}=(-1)^{n-\lambda_1}
f(\lambda)$ holds 
for any  $\lambda=(\lambda_1,\dots,\lambda_r)
\vdash n$.
For the case $\lambda_i\ge 2$,
by taking $\mu=\lambda'$ 
and $\mu\upo{i}=\lambda$,
Theorem~\ref{the4} 
follows from Proposition~\ref{pp4}
and the fact $\eta_{\lambda}=(-1)^{n-\lambda_1}
f(\lambda)$
directly. 

Now consider the case
$\lambda_i=1$. 
Then due to the assumption of $\lambda$, we have $i=s$. Thus $\lambda=(\lambda_1,\dots,\lambda_{i-1},1)$ and $\lambda'=(\lambda_1,\dots,\lambda_{i-1})$. As $i=s\ge 2$, by (\ref{def1}),
\aln{tha-e1}
{
	f(\lambda)-f(\lambda')=f(\lambda-\hat{1}).
}
Hence the result follows from 
the fact  $\eta_{\lambda}=(-1)^{n-\lambda_1}
f(\lambda)$. 
\proofend

\section{Proof of Theorem~\ref{the6}
	\label{secthree}}

In this section, we further develop Proposition~\ref{pp4} to prove Theorem~\ref{the6}.

	For $r\ge 2$ and $2\le i<j\le r$, let 
	$\P_n(r,i,j)$ denote the set of partitions $\lambda=(\lambda_1,\dots,\lambda_r)
	$ in $\P_n(r,i)$ with $\lambda_j>\lambda_{j+1}$ if $j<r$.
    For any $\lambda=(\lambda_1,\dots,\lambda_r)
	\in \P_n(r,i,j)$, let $\lambda(i,j)$ denote the partition
	$(\lambda \upo{i})\downo{j}$.
	Clearly, $\lambda(i,j)\vdash n$, and $\ell(\lambda(i,j))\le \ell (\lambda)$, where the inequality is strict if and only if $j=r$ and $\lambda_j=1$.
	
In what follows, we focus on the difference between $f(\mu)$ and $f(\mu(i,j))$ for any partition $\mu\in\P_n(s,i,j)$. We first give a recurrence relation for the case when $j=s$.

\begin{lemm}\label{pp5}
For any $\mu=(\mu_1,\dots,\mu_s)\in \P_n(s,i,s)$, where $2\le i\le s-1$, 
\eqn{in1-8}
{
f(\mu(i,s))-f(\mu)
	 &=&(2\mu_i-2\mu_s+s-i+2) 
	 f(\mu\upo{i}-\hat{1})
	 \nonumber \\
	& & -(2\mu_i+s-i)f(\mu-\hat{1})
	 +2(\mu_s-1)f(\mu(i,s)-\hat{1}).
}
\end{lemm}

\proof 
When $\mu_s=1$, (\ref{def1}) implies that
\aln{pp5-e0}
{
f(\mu\upo{i})=f(\mu\upo{i}\setminus \mu_s)+f(\mu\upo{i}\setminus \mu_s-\hat{1})=f(\mu(i,s))+f(\mu\upo{i}-\hat{1}).
}
Then applying Proposition~\ref{pp4}, we have
\aln{pp5-e3}
{
f(\mu\upo{i})-f(\mu)=(2\mu_i+s-i+1)f(\mu\upo{i}-\hat{1})-(2\mu_i+s-i)f(\mu-\hat{1}).
}
Thus (\ref{pp5-e0}) and (\ref{pp5-e3}) together give us
\aln{pp5-e4}
{
f(\mu(i,s))-f(\mu)=(2\mu_i+s-i)f(\mu\upo{i}-\hat{1})-(2\mu_i+s-i)f(\mu-\hat{1}).
}
Hence the result holds when $\mu_s=1$.

Now we assume $\mu_s\ge 2$. Then $\mu,\mu\upo{i},\mu(i,s)$ have the same length $s$.
Again, by applying Proposition~\ref{pp4}, we have 
\equ{pp5-e1}
{
f(\mu\upo{i})-f(\mu(i,s))
=(2(\mu_s-1)+1)f(\mu\upo{i}-\hat{1})
-2(\mu_s-1)f(\mu(i,s)-\hat{1})
}
and 
\equ{pp5-e2}
{
	f(\mu\upo{i})-f(\mu)
	=(2\mu_i+s-i+1)f(\mu\upo{i}-\hat{1})
	-(2\mu_i+s-i)f(\mu-\hat{1}).
}
Thus, (\ref{in1-8}) follows from (\ref{pp5-e1}) and (\ref{pp5-e2}).
\proofend 

Clearly, $\mu\triangleleft \mu(i,j)$ holds for any $\mu\in\P_n(s,i,j)$.
By applying Lemma~\ref{pp5}, we can show that 
$f(\mu(i,j))=f(\mu)$ holds whenever 
$\mu\in \P^*_{n,3}$ (i.e., 
	$\mu_1=3$ and $\mu_2\le 2$)
and the $i$-th part of $\mu$ is $1$.

\begin{coro}\label{co3}
For $\mu=(\mu_1,\dots,\mu_s)\in \P_n(s,i,j)$, if 
$\mu\in \P^*_{n,3}$ and $\mu_i=1$,
then $f(\mu(i,j))=f(\mu)$. 
\end{coro}

\proof 
The given conditions on $\mu$ imply that $\mu_t=2$ for $2\le t\le i-1$
	and 
	$\mu_t=1$ for $i\le t\le s$.
Since $\mu\in \P_n(s,i,j)$, 
either $j=s$ or $\mu_j>\mu_{j+1}$.
Thus, $j=s$.
As $\mu_i=\mu_s=1$, 
by Lemma~\ref{pp5}, we have 
\equ{co3-e1}
{
f(\mu(i,j))-f(\mu)
=f(\mu(i,s))-f(\mu)
=(s-i+2) 
\big ( f(\mu\upo{i}-\hat{1})
-f(\mu-\hat{1})\big ).
} 
By (\ref{def1}), 
\equ{co3-e2}
{
f(\mu\upo{i}-\hat{1})
-f(\mu-\hat{1})
=f(\underbrace{2,1,\dots,1}_{i-1},1)
-f(\underbrace{2,1,\dots,1}_{i-1})
=f(1)=0.
}
Hence the conclusion follows.
\proofend

\textbf{Remark.}
It can be proved by induction that 
for any $\mu\in\P^*_{n,3}$, 
$f(\mu)
=2n+2$ holds.

On the other hand, applying Proposition~\ref{pp4}, 
we can conclude 
that $f(\mu\upo{i})>f(\mu)$ for any partition $\mu\in\P_n(s,i)$.
It will be applied to show that 
$f(\mu\upo{i}-\hat{1})\ge f(\mu-\hat{1})$.

\begin{lemm}\label{lem4-1}
	Let $\mu=(\mu_1,\dots,\mu_s)\in\P_n(s,i)$, where $2\le i\le s$. Then
	\aln{lem4-1-e1}
	{
		f(\mu\upo{i})-f(\mu)\ge 
		f(\mu\upo{i}-\hat{1})>0.
	}
\end{lemm}

\proof Obviously, we need only to consider $n\ge 3$. 
If $n=3$, then $\mu=(2,1)$, $s=i=2$, 
$\mu\upo{i}=(2,2)$, implying that 
$$
f(\mu\upo{i})-f(\mu)=f(2,2)-f(2,1)=5-2
>1=f(1,1)=f(\mu\upo{i}-\hat{1}).
$$
Assume that $n\ge 4$ and 
the conclusion holds 
for all partitions $\mu\in \P_{n'}(s',i')$,
where $n'\le n-1$ and $2\le i'\le s'$.
Now let $\mu\in \P_{n}(s,i)$
and $2\le i\le s$.
By Proposition~\ref{pp4}, 
\equ{lem4-1-e2}
{
	f(\mu\upo{i})-f(\mu)
	=f(\mu\upo{i}-\hat{1})+
	(2\mu_i+s-i)
	\big (f(\mu\upo{i}-\hat{1})
	 -
	f(\mu-\hat{1})\big ).
}
As $\ell(\mu\upo{i}-\hat{1}) \ge 2$, 
we have $f(\mu\upo{i}-\hat{1})>0$.

If $\mu_i\ge 2$, then 
$\mu\upo{i}-\hat{1}=(\mu-\hat{1})\upo{i}$.
By induction, 
$f(\mu\upo{i}-\hat{1})>f(\mu-\hat{1})$ holds for this case.
Now assume that $\mu_i=1$.
Then, 
$$
\mu\upo{i}-\hat{1}
=(\mu_1-1,\dots,\mu_{i-1}-1,1),
\quad 
\mu-\hat{1}
=(\mu_1-1,\dots,\mu_{i-1}-1).
$$
By (\ref{def1}), 
\equ{lem4-e3}
{
f(\mu\upo{i}-\hat{1})-
f(\mu-\hat{1})
=f((\mu-\hat{1})-\hat{1})\ge 0.
}
The conclusion holds for $\mu$.
\proofend 

\begin{lemm}\label{lem4}
Let $\mu=(\mu_1,\dots,\mu_{s})\in\P_n(s,i)$, where $2\le i\le s$.
	Then 
	$f(\mu\upo{i}-\hat{1})\ge f(\mu-\hat{1})$, where the equality holds if and only if 
$\mu\in \P^*_{n,3}$ and $\mu_i=1$.
\end{lemm}

\proof
If $\mu\upo{i}-\hat{1}$ and $\mu-\hat{1}$ have the same length, then $\mu\upo{i}-\hat{1}=(\mu-\hat{1})\upo{i}$, and thus $f(\mu\upo{i}-\hat{1})> f(\mu-\hat{1})$ follows from Lemma~\ref{lem4-1}. 
In the following, we consider the case when $\ell(\mu\upo{i}-\hat{1})>\ell(\mu-\hat{1})$, or equivalently, $\mu_i=1$. 

Since $\mu\in\P_n(s,i)$, we have $\mu_{i-1}>\mu_i$. Thus $\mu\upo{i}-\hat{1}=(\mu_1-1,\dots,\mu_{i-1}-1,1)$ and $\mu-\hat{1}=(\mu_1-1,\dots,\mu_{i-1}-1)$. Since $i\ge 2$, by (\ref{def1}),
\eqn{tha-e7}
{
f(\mu\upo{i}-\hat{1})-f(\mu-\hat{1})=
f((\mu-\hat{1})-\hat{1})
\ge 0,
}
where the equality holds if and only if $(\mu-\hat{1})-\hat{1}=(1)$,
i.e., 
$$
\mu_1=3, \mu_2=\cdots=\mu_{i-1}=2,
 \mbox{ and }\mu_i=\cdots=\mu_s=1.
 $$
Hence the result holds.
\proofend

Now we compare the value of $f(\mu(i,j))$ with $f(\mu)$ for any $\mu\in \P_n(s,i,j)$.

\begin{prop}\label{th1-1}
Let $\mu=(\mu_1,\dots,\mu_s)\in \P_n(s,i,j)$, where $2\le i<j\le s$.
Then
	\eqn{in1-8-0}
	{
		f(\mu(i,j))\ge f(\mu),
	}
where the equality holds if and only if
$\mu\in \P^*_{n,3}$ and $\mu_i=1$.
\end{prop}

\proof 
By Corollary~\ref{co3}, 
we need only to prove that 
$f(\mu(i,j))> f(\mu)$ whenever $\mu_1\ne 3$, or $\mu_2\ge 3$ 
or $\mu_i\ge 2$.

Now suppose the conclusion fails, 
and 
$n$ is the minimum integer in $\N$ 
with some $\mu=(\mu_1,\dots,\mu_s)\in \P_n(s,i,j)$, where 
$\mu_1\ne 3$, or $\mu_2\ge 3$ 
or $\mu_i\ge 2$,
such that 
$f(\mu(i,j))\le f(\mu)$. 
We will complete the proof by showing the following claims. 

\inclaim $\mu_j\ge 2$.

Suppose the claim fails, i.e., $\mu_j=1$. Then $j=s$, as $\mu\in\P_n(s,i,j)$ implies that either 
$j=s$ or $\mu_{j}>\mu_{j+1}\ge 1$.
By Lemma~\ref{pp5}, 
\aln{th1-e1}
{
	f(\mu(i,j))-f(\mu)=(2\mu_i+s-i) 
	\big(f(\mu\upo{i}-\hat{1}) - f(\mu-\hat{1})\big ).
}

Then by Lemma~\ref{lem4} and the 
assumption on $\mu$, 
 (\ref{th1-e1}) implies that 
 $f(\mu(i,j))>f(\mu)$, a contradiction.
\claimend

Claim 1 implies that $\mu(i,j)$ and $\mu$ have the same length.

\inclaim $j<s$.

Suppose that $j=s$. Then, by 
Lemma~\ref{pp5}, 
\eqn{th1-1-e1}
{
f(\mu(i,s))-f(\mu)
	&=&(2\mu_i-2\mu_s+s-i+2) 
\big(f(\mu\upo{i}-\hat{1}) - f(\mu-\hat{1})\big )
\nonumber \\
& &
	+2(\mu_s-1)
\big(	f(\mu(i,s)-\hat{1})-f(\mu-\hat{1})\big ).
}
By Lemma~\ref{lem4} and the assumption on $\mu$, 
$f(\mu\upo{i}-\hat{1}) >f(\mu-\hat{1})$.

Note that $j=s$.
Claim 1 implies $\mu_s\ge 2$.
Thus,
$\mu(i,s)-\hat{1}=(\mu-\hat{1})(i,s)$.
By the assumption on the minimality of $n$, 
$f(\mu(i,s)-\hat{1})\ge f(\mu-\hat{1})$
holds. 

Thus, by (\ref{th1-1-e1}), 
the above conclusions imply that 
$f(\mu(i,s))>f(\mu)$, 
a contradiction to the assumption. 
\claimend

\inclaim $f(\mu(i,j))-f(\mu)> 0$.

By Claim 2 and the assumption that $\mu\in \P_n(s,i,j)$, $\mu_j>\mu_s$. Then for any $0\le k\le \mu_s$, $\mu(i,j) \setminus \mu_s-\hat{k}$ and $\mu\setminus \mu_s-\hat{k}$ have the same size, which implies that 
$\mu(i,j)\setminus \mu_s-\hat{k}=(\mu\setminus \mu_s-\hat{k})(i,j)$.

By (\ref{def1}),  we have 
\eqn{th1-1-e2}
{
	f(\mu(i,j))-f(\mu)
	&=&f(\mu(i,j)\setminus \mu_s)
	-f(\mu\setminus \mu_s)
	\nonumber \\
	& &
+\sum_{k=1}^{\mu_s}
	 \binom{\mu_s}{k}(2k-1)!!
	\big ( f(\mu(i,j) \setminus \mu_s-\hat{k})
	-f(\mu\setminus \mu_s-\hat{k}) \big)
	\nonumber\\
&=&f((\mu\setminus \mu_s)(i,j))
-f(\mu\setminus \mu_s)
\nonumber \\
& &
+\sum_{k=1}^{\mu_s}
\binom{\mu_s}{k}(2k-1)!!
\big ( f((\mu\setminus \mu_s-\hat{k})(i,j) )
-f(\mu\setminus \mu_s-\hat{k}) \big).
}

By the assumption on the minimality of $n$, the following inequalities hold, among which at most one equality holds:
$$
f((\mu\setminus \mu_s)(i,j))\ge f(\mu\setminus \mu_s),
\qquad 
f((\mu\setminus \mu_s-\hat{k})(i,j))
\ge f(\mu\setminus \mu_s-\hat{k}),
\quad \forall k:1\le k\le \mu_s.
$$
By (\ref{th1-1-e2}), Claim 3 holds,
contradicting the assumption of $\mu$. 
\claimend

Hence Proposition~\ref{th1-1} follows.
\proofend

Now we conclude this section by a proof of Theorem~\ref{the6}.

\noindent\textit{Proof of Theorem~\ref{the6}.}
We need only to consider the case when $2\le u< n$.
	For any $\lambda,\lambda'\in 
	\P_{n,u}$ with $\lambda \triangleleft \lambda'$, 
	there is a sequence of partitions
	$\gamma^1(=\lambda),\gamma^2, \dots,\gamma^t(=\lambda')$ in $\P_{n,u}$ 
	such that for any $q$ with $1\le q<t$, 
	$\gamma^{q+1}=\gamma^{q}(i,j)$ holds 
	for some $i,j$ with $2\le i<j\le r$, where 
	$r=\ell(\gamma^q)$.
	Also see \cite{ku2013} for this conclusion.
Then the result follows from 
Proposition~\ref{th1-1}.
\proofend

\section{Further Study}

For any $u_1,u_2\in\N$ with $u_1\ge u_2$
and $a,b\ge 0$, 
denote by $(u_1^a,u_2^{b})$ the partition 
	$\mu=(\mu_1,\dots,\mu_{a+b})$,
	where 
$\mu_t=u_1$ for $1\le t\le a$
and  $\mu_t=u_2$
for $a+1\le t\le a+b$.
We omit the $a$ (or resp., $b$) if $a=1$ (or resp., $b=1$).
It can be proved by induction that 
for $a\ge 1$ and $b\ge 0$, 
\equ{}
{
f(2^a,1^b)=a^2+b(a-1)+1.
}
It is known that 
$f(\mu)=2n+2$ for each $\mu\in \P^*_{n,3}$.
Then, it is not difficult to verify that 
for any $n\ge 10$, if $4\le a\le \frac n2$,
then $f(\mu)<f(2^a,1^{n-2a})$ holds for 
each $\mu\in \P^*_{n,3}$,
although $(2^a,1^{n-2a})\triangleleft \mu$	whenever $\mu$ has at most $n-2a-1$ parts equal to $1$.

More general, it can be shown by induction (in a similar manner to the proof of Proposition~\ref{pp4}) that for any $u,q\in\N$,
the following two identities hold:
\eqn{}
{
f(u+1,u^{q-1})
&=&2u(f(u,(u-1)^{q-1})+f((u-1)^q)),
\label{eq5-1}
\\
f((u+2)^q,1)&=&(2u+q+1)f((u+1)^{q},1)+2uf(u^{q},1).
\label{eq5-2}
}
Applying (\ref{eq5-1}) and (\ref{eq5-2}),
it can be further proved by induction that
\equ{}
{
qf(u+1,u^{q-1})=2u f(u^q,1).
}
As a result, $f(u^q,1)>f(u+1,u^{q-1})$ when $q>2u$, while 
$(u^q,1)\triangleleft(u+1,u^{q-1})$.

These above observations imply that the condition $\lambda,\lambda'\in 
\P_{n,u}$ for some $u\in\N$ in Theorem~\ref{the6} is crucial. However, we wonder whether the following conjecture holds.

\begin{conj}
	For $n\ge 2$, $\lambda\in 
	\P_{n,u}$ and $\mu\in \P_{n,v}$ for some $u,v\ge 2$ with $v\ge u+2$, 
	if $\lambda \triangleleft \mu$,
	then
	$|\eta_{\lambda}|	< |\eta_{\mu}|$.
\end{conj}

\section*{Acknowledgement}

This research is supported by the Ministry of Education,
Singapore, under its Academic Research Tier 1 (RG19/22). Any opinions,
findings and conclusions or recommendations expressed in this
material are those of the authors and do not reflect the views of the
Ministry of Education, Singapore.
The first author would like to express her gratitude to National Institute of Education and Nanyang Technological University of Singapore for offering her Nanyang Technological University Research Scholarship during her PhD study.

\noindent {\bf Data availability}\hspace{0.2 cm} Data sharing not applicable to this article as no datasets were generated or analysed during the current study.

\end{document}